\documentclass{elsart}

\usepackage{amsmath,amssymb,graphicx,psfrag}

\newcommand{\deadlamp}{deadlamp}
\newcommand{\seesaw}{seesaw}
\newcommand{\lamptat}{lamptat}

\newcommand{\ctcomb}{ctcomb}
\newcommand{\ct}{ct}
\newcommand{\Epstein}{MR93i:20036}
\newcommand{\Gilman}{MR97a:20054}
\newcommand{\Grig}{MR2002j:60009}
\newcommand{\GrigNag}{GrigNag}

\newcommand{\gubagrowth}{gubagrowth}
\newcommand{\NSnotes}{NSnotes}
\newcommand{\NSpaper}{MR1329042}
\newcommand{\Sipser}{Sipser}
\newcommand{\ElstonOstheimer}{EO}
\newcommand{\MeBSCF}{MeBSCF}
\newcommand{\belkbrown}{belkbrown}
\newcommand{\josegrowth}{josegrowth}
\newcommand{\cfp}{cfp}
\newcommand{\blakegd}{blakegd}
\newcommand{\cannoncone}{cannoncone}
\newcommand{\belkbux}{belkbux}
\newcommand{\MS}{mitrana}

\newcommand{\bi}{\begin{itemize}}
\newcommand{\ei}{\end{itemize}}
\newcommand{\be}{\begin{enumerate}}
\newcommand{\ee}{\end{enumerate}}
\newcommand{\bc}{\begin{center}}
\newcommand{\ec}{\end{center}}
\newcommand{\bt}{\begin{tabular}}
\newcommand{\et}{\end{tabular}}

\newcommand{\gset}{generating set}
\newcommand{\cf}{context-free}
\newcommand{\lan}{language}

\newcommand{\pda}{pushdown automaton}

\newcommand\Z{\mathbb Z}

\begin{document}

\begin{frontmatter}

\title{Cone types and geodesic
languages for\\ lamplighter groups and Thompson's group $F$}

\author{Sean Cleary\thanksref{sean}} 

\address{Department of
Mathematics, The City College of New York, New York} 

\ead{cleary@sci.ccny.cuny.edu}

\thanks[sean]{The first author acknowledges support from PSC-CUNY grant 
\#65752-0034}

\author{Murray Elder\thanksref{murray}}

 \address{School of Mathematics and
Statistics, University of St Andrews, Scotland} 

\ead{murray@mcs.st-and.ac.uk}

\thanks[murray]{The second
author acknowledges support from EPSRC grant GR/S53503/01}

\author{Jennifer Taback\thanksref{jen}}

\address{Department of
Mathematics, Bowdoin College, Brunswick, Maine} 

\ead{jtaback@bowdoin.edu}

\thanks[jen]{The third author
acknowledges support from NSF grant DMS-0437481}

\begin{abstract}
We study languages of geodesics in lamplighter groups and Thompson's
group $F$.  We show that the lamplighter groups $L_n$ have infinitely
many cone types, have no regular geodesic languages, and have 1-counter,
context-free and counter geodesic languages with respect to certain generating
sets.  We show that the full language of geodesics with respect to
one generating set for the lamplighter group is not counter but is context-free,
while with respect to another generating set the full language of geodesics is counter
and context-free.  In Thompson's group $F$  with respect to the standard finite
generating set, we show there are infinitely many cone types and that there is
no regular language of geodesics.  We show that the existence of families of ``seesaw'' elements
with respect to a given generating set in a finitely generated
infinite group precludes a regular language of geodesics and
guarantees infinitely many cone types with respect to that generating
set.
\end{abstract}

\end{frontmatter}



\section{Introduction}\label{sec:intro}

In this article we prove some language-theoretic consequences of
work by Cleary and Taback  describing geodesic words in the
lamplighter groups  \cite{\deadlamp,\lamptat} and Thompson's group $F$
\cite{\ctcomb,\seesaw}.  We consider Thompson's group $F$ with its standard
finite generating set,
$$F = \langle x_0,x_1
|[x_0x_1^{-1},x_0^{-1}x_1x_0],[x_0x_1^{-1},x_0^{-2}x_1x_0^2]
\rangle.$$  We consider the lamplighter groups $L_m$ in the
standard wreath product presentation,
$$
L_m = \langle a,t | [t^{i}a t^{-i},t^{j}at^{-j}], a^m, i,j \in \Z
\rangle.$$  We also consider an alternate generating set for $L_2$
which arises from considering this group as an automata group, in
which case the natural generators are $t$ and $(ta)$.

We prove the following results for $L_m$ with $m \geq 2$ with respect
to generating set $\{a,t\}$, for $L_2$ with respect to the automata
group generating set $\{t,(ta)\}$, and for Thompson's group $F$
with respect to the standard finite generating set $\{x_0, x_1\}$.

\bi
\item There are infinitely many cone types; that is, there are infinitely
many families of possible geodesic extensions to group elements.
\item There is no regular language of geodesics which includes at least
one representative of each group element.
\ei

We prove that there are 1-counter languages of geodesics with a unique
representative for each element for $L_m$ with $m \geq 2$ with respect
to generating set $\{a,t\}$ and for $L_2$ with respect to the automata
group generating set $\{t,(ta)\}$.

However, the formal language class of the full language of geodesics
in $L_2$ depends on the choice of \gset. For the automaton \gset\
$\{t,ta\}$ the set of all geodesics is 1-counter, which implies it
is both context-free and counter, yet for the wreath generating set
$\{a,t\}$ the set of all geodesics for $L_m$ with $m \geq 2$ is shown
to be context-free but not counter.

In addition, we show that if a finitely generated group contains a
family of ``seesaw elements'' of arbitrary swing with respect to a
finite generating set, then the group has infinitely many cone types
with respect to that generating set
and cannot have a regular language of geodesics with respect to that generating
set.

The concepts of cone types, regular, \cf, counter and 1-counter
geodesic languages are intimately connected, and tell us much about
the structure of geodesics in a given group presentation. If a group
has finitely many cone types, then the full \lan\ of geodesics is
regular. The converse is true when all relators in the presentation
have even length \cite{\NSnotes}, and is conjectured to be true in
general.  In addition, Grigorchuk and Smirnova-Nagnibeda \cite{\GrigNag}
showed that if a group has finitely many cone types then it
has a rational complete growth function. While the
growth of the lamplighter groups is much studied, the rationality of
the growth function for Thompson's group is an open question.

We would like to thank Katherine St.~John, Gretchen Ostheimer and the 
anonymous referee
for helpful comments on drafts of this article.

\section{Languages, Grammars and Cone types}\label{sec:defn}

We now present the necessary definitions for the different types
of languages and finite state automata which we consider below. If
$X$ is a finite set of symbols then we let $X^*$ be the set of all
finite strings in the symbols of $X$, including the {\em empty
string} $\epsilon$ which has no letters.  A language over $X$ is a
subset of $X^*$.  A regular language is one accepted by some
finite state automaton. See Epstein {\it et al} \cite{\Epstein} and
Sipser \cite{\Sipser} for an
introduction to regular languages.

A key tool in the theory of regular \lan s is the Pumping Lemma:

\begin{lem}[Pumping Lemma for regular \lan s]\label{lem:pumpR}
Let $A$ be a regular \lan. There is a number $p$ (the pumping length)
so that if $s$ is any word in $A$ of length at least $p$, then $s$ may
be divided into three pieces $s=xyz$ where we can require either
that $|xy|\leq p$ or $|yz|\leq p$ and have both \bi
\item For each $i\geq 0$, $xy^iz\in A$ \item $|y|>0$ \ei
\end{lem}
See Sipser \cite{\Sipser} for a proof of the Pumping Lemma. A standard
method of showing that a given \lan\ fails to be regular is via
this lemma.

Let $G$ be a group with a word metric defined with respect to a
finite generating set.  Cannon \cite{\cannoncone} defined the cone type
of an element $w \in G$ to be the set of geodesic extensions of $w$ in
the Cayley graph.
\begin{defn}[Cone type]
A path $p$ is {\em outbound} if $d(1,p(t))$ is a strictly
increasing function of $t$. For a given $g\in G$, the {\em cone}
at $g$, denoted $C'(g)$ is the set of all outbound paths starting
at $g$. The {\em cone type} of $g$, denoted $C(g)$, is
$g^{-1}C'(g)$.
\end{defn}

This definition applies both in the discrete setting of the group
and in the one-dimensional metric space which is the Cayley graph.
In the Cayley graph,  a cone type may include paths that
end at the middle of an edge. If the presentation of $G$ consists
of only even-length relators, then all cone types will consist
solely of full edge paths. Neumann and Shapiro \cite{\NSpaper} show that
whether a group has finite or infinitely many cone
 types may depend on the choice of generating set  and
 they describe a number of properties in of cone types in \cite{\NSnotes}.

In all that follows, when we consider the set of cone types of a
group, we fix a particular generating set.  Knowing the set of all
possible cone types of elements of a particular group provides
information about the full language of geodesics of that group
with respect to that particular generating set. We may also get
information and about the growth function of the group as follows.

\begin{lem}\cite{\NSnotes}\label{lem:coneimpliesreg}
If $G$ has finitely many cone types with respect to a finite
generating set, then the full language of geodesics with respect to
that \gset\ is regular.
\end{lem}

It is not known whether the converse is true in general, but the
converse is true for presentations having only even-length
relators.

\begin{lem} \cite{\NSnotes}\label{lem:evenlength}
Let $G$ be a finitely-presented group with a generating set in which all
relators have even length.  If the full language of geodesics with
respect to this generating set is regular, then $G$ has finitely
many cone types with respect to this generating set.
\end{lem}

We consider not only regular languages below, but context
free and counter languages as well.

\begin{defn}[Pushdown automaton]
A {\em pushdown automaton} is a machine consisting of
\bi
\item a stack,
\item a finite set of stack symbols, including
$\$$ which is a marker symbol for the beginning and end of the
stack, and 
 \item a finite state automaton with possible additional edge labels to
 affect the stack \ei
where each transition may, in addition to reading a letter from the input word,
potentially also push or pop a stack symbol on or off
the stack. A word is accepted by the pushdown automaton if it
represents a sequence of transitions from the start state to an
accept state so that the stack is empty at the final state.
\end{defn}

\begin{defn}[Context-free]
A language is {\em context-free} if it is the set of all strings
recognized by some \pda.
\end{defn}

\begin{figure}[ht!]
\bc
\psfrag{e}{$\epsilon$}
\psfrag{q0}{$q_0$}
\psfrag{A}{$A$}
\psfrag{e,push s}{$(\epsilon$,push $\$)$}
\psfrag{e,pop s}{$(\epsilon$,pop $\$)$}
\psfrag{a,push 1}{$(a$,push $1)$}
\psfrag{b,pop 1}{$(b$,pop $1)$}

  \includegraphics[height=3cm]{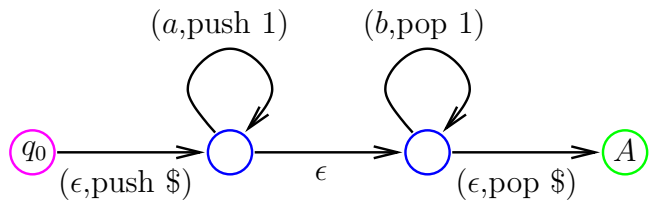}
  \ec
  \caption{Pushdown automaton accepting $a^nb^n$ with start state $q_0$.}
  \label{fig:PDAanbn}
\end{figure}

For example, the language $\{a^nb^n \; | \; n\in \mathbb N\}$ is
accepted by the \pda\ in Figure \ref{fig:PDAanbn} with
alphabet $a,b$ and stack symbols $\$,1$, and this language is not
regular.

The following result is proved in \cite{\Sipser}.
\begin{lem}\label{lem:CFunion}
Context-free languages are closed under union.
\end{lem}

There is a pumping lemma for context-free languages as well.

\begin{lem}[Pumping Lemma for \cf\ \lan s]\label{lem:pumpCF}
Let $A$ be a \cf\ \lan. There is a number $p$ (the pumping
length) so that if $s$ is any word in $A$ of length at least $p$,
then $s$ may be divided into five pieces $s=uvxyz$ such that \bi
\item For each $i\geq 0$, $uv^ixy^iz\in A$ \item $|v|,|y|>0$ \item
$|vxy|\leq p$ \ei
\end{lem}
See \cite{\Sipser} for a proof of Lemma \ref{lem:pumpCF}.

\begin{defn}[$G$-automaton]\label{defn:EE} 
Let $G$ be a group and $\Sigma$ a finite set.  A (non-deterministic)
 {\em $G$-automaton} $A_G$ over $\Sigma$ is a finite directed graph
 with a distinguished {\em start vertex} $q_0$, some distinguished
 {\em accept vertices}, and with edges labelled by elements of $(\Sigma^{\pm
 1}\cup\{\epsilon\})\times G$.  If $p$ is a path in $A_G$, then the element
 of $(\Sigma^{\pm 1})$ which is the first component of the label of
 $p$ is denoted by $w(p)$, and the element of $G$ which is the second
 component of the label of $p$ is denoted $g(p)$. If $p$ is the empty
 path, $g(p)$ is the identity element of $G$ and $w(p)$ is the empty
 word.  $A_G$ is said to {\em accept} a word $w\in (\Sigma^{\pm 1})$
 if there is a path $p$ from the start vertex to some accept vertex
 such that $w(p)=w$ and $g(p)=_G 1$.
 \end{defn}

If $G$ is the trivial group, then a $G$-automaton is just a finite
state automaton.
\begin{defn}[Counter]
A language is {\em $k$-counter} if it is accepted by some
$\Z^k$-automaton. We call the (standard) generators of $\Z^k$ and
their inverses {\em counters}. A language is {\em counter} if it is
$k$-counter for some $k\geq 1$.
\end{defn}

\begin{figure}[ht!]
  \begin{center}
  \psfrag{e}{$\epsilon$}
\psfrag{q0}{$q_0$}
\psfrag{A}{$A$}
\psfrag{a,xy}{$(a,x_1x_2)$}
\psfrag{b,X}{$(b, x_1^{-1})$}
\psfrag{a,Y}{$(a, x_2^{-1})$}

   \includegraphics[height=3cm]{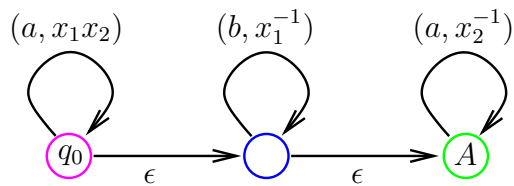}
  \end{center}
  \caption{A counter automaton accepting $a^nb^na^n$.}
  \label{fig:counternotCF}
\end{figure}

For example, the language $\{a^nb^na^n \; | \; n\in \mathbb N\}$ is
accepted by the $\Z^2$-automaton in Figure \ref{fig:counternotCF},
with alphabet $a,b$ and counters $x_1,x_2$, and this language is not
\cf\ \cite{\Sipser}.

In the case of $\Z$-automata, we assume that the generator is $1$ and
the binary operation is addition, and we may insist without loss of
generality that each transition changes the counter by either $0,1$ or
$-1$. We can do this by adding states and transitions to the automaton
appropriately.  That is, if some edge changes the counter by $k\neq
0,\pm 1$ then divide the edge into $|k|$ edges using more states.  The
symbols $+,-$ indicate a change of $1,-1$ respectively on a
transition.

Note that with this definition, a $\Z$-automaton cannot ``see'' the
value of the counter until it reaches an accept state, so is not
equivalent to a \pda\ with one stack symbol, which can determine if
the stack is empty at any time.  In the literature these counter
automata are sometimes called ``partially blind'' to emphasize this
difference.

Elder \cite {\MeBSCF} shows that counter languages have the following
properties:
 
\begin{lem}
1-counter languages are context-free.
\end{lem}

\begin{lem}[Closure properties]
\label{lem:closurekcounter} 
1-counter languages are not closed under concatenation or
intersection.  If $C$ is $k$-counter for $k\geq 1$ and $L$ is regular,
then $C\cap L$, $CL$ and $LC$ are all $k$-counter. The union of a
finite number of $k$-counter languages is $k$-counter.
\end{lem}

 In order to show that a language is not counter we make use of the
following lemma from  \cite {\MeBSCF}.
\begin{lem}[Swapping Lemma] \label{lem:swap}
If $L$ is counter then there is a constant $s>0$, the ``swapping
length'', such that if $w\in L$ with length at least $2s+1$ then $w$
can be divided into four pieces $w=uxyz$ such that $|uxy|\leq 2s+1$,
$|x|,|y|>0$ and $uyxz\in L$.
\end{lem}

\textit{Proof}
Let $p$ be the number of states in the counter automaton. If a
path visits each state at most twice then it cannot have length
more than $2p$, so $w$ visits some state $s$ at least three times.
Let $u$ be the prefix of $w$ ending when $w$ reaches state $s$ for
the first time.  Let $x$ be the continuation of $u$, ending when
$w$ reaches state $s$ for the second time.  Let $y$ be the
continuation of $x$, ending when $w$ reaches state $s$ for the
third time.  Finally, let $z$ be the remainder of $w$.  So
$w=uxyz$ ends at an accept state, with all the counters balanced
correctly. If we switch the orders of the strings $x$ and $y$, then we
will also have a path leading to the same accept state with the
correct counters, so $uyxz\in L$.
\hfill$\Box$

The Swapping Lemma is similar to the Pumping Lemmas for regular
and context-free languages.  It is only useful if the word $w$ has no
repeated concurrent subwords; if there are repeated subwords, we can
take $x=y$ in the statement of the lemma, and the resulting word will
be identical to the initial word so the conclusion is vacuously true.

For more background on counter \lan s see
Mitrana and Stiebe \cite{\MS}, Elston and Ostheimer
\cite{\ElstonOstheimer},  Elder \cite{\MeBSCF}, and Gilman
\cite{\Gilman}.

\section{Geodesic languages for the lamplighter groups in 
the wreath product \gset}
\label{sec:lampwreath}

The lamplighter group $L_m$, with presentation
$$L_m = \langle a,t | [t^{i}a t^{-i},t^{j}at^{-j}], a^m \rangle$$ is
the wreath product $\Z_m \wr \Z$, where the generator $a$ in the
above presentation generates $\Z_m$ and $t$ generates $\Z$. We
will refer to these generators as the wreath product generators to
distinguish them from the generators which arise naturally when
the group $L_2$ is considered as an automata group.

In \cite{\deadlamp}, Cleary and Taback studied metric properties of
this wreath product presentation for $L_m$. In this paper we
consider some language theoretic consequences of that work.  We
begin with a geometric interpretation for elements of $L_m$, first
in the case $m=2$, and then for $m >2$.

An element of $L_2$ is best understood via the following
geometric picture.  We visualize each $\Z_2$ factor as a light bulb which
is either on or off, and the wreath product with $\Z$ thus creates an
infinite string of light bulbs, one at each integer.  An element
$w \in L_2$ is represented by a configuration of a finite number
of illuminated bulbs, and a position of the ``lamplighter'' or
cursor.  A word $\gamma$ in the generators $\{a,t\}$ denoting the
element $w$ can be thought of as a sequence of instructions for
creating the configuration representing $w$, in the following way.

The position of the cursor indicates the current bulb under
consideration by the ``lamplighter''.  The generator $a$ changes the
state of the bulb at the cursor position, and the generator $t$
moves the cursor one unit to the right.  Thus prefixes of the word
$\gamma$ appear to ``move'' the cursor to different integral
positions, possibly changing the state of bulbs as well, ending
with the configuration representing $w$.

\begin{figure}
\bc
\includegraphics[width=13cm]{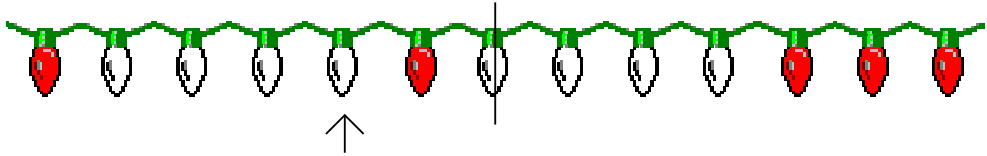}
\ec
\caption{The element $w=a_4 a_5 a_6 a_{-1} a_{-6} t^{-2}$ of $L$,
expressed as a configuration of illuminated bulbs and a cursor in
the position $-2$.  Shaded bulbs represent bulbs which are
illuminated, clear bulbs represent bulbs which are off, the
vertical bar denotes the origin in $\Z$, and the arrow denotes the
position of the cursor. \label{fig:lampex1}}
\end{figure}

The identity word is represented by the configuration of bulbs
which are all in the off state, and the cursor at the origin.
Figure \ref{fig:lampex1} gives an example of an element of $L_2$
represented in this way.  For a given word $w$, we consider the
successive prefixes of $w$ as steps involved in the creation of
$w$.

Elements of $L_m$ for $m>2$ can be understood via the analogous
picture in which the bulbs have $m$ states.  Again, the generator
$t$ moves the cursor one unit to the right, and $a$ increments the
state of the current bulb under consideration.

\subsection{Normal forms for elements of $L_2$}

We first give normal forms for elements of $L_2$ in detail, then
generalize these forms to $L_m$ for $m > 2$.

Normal forms for elements of $L_2$ with respect to the wreath
product generators $a$ and $t$ are given in terms of the
conjugates $a_k=t^k a t^{-k}$, which move the cursor to the $k$-th
bulb, turn it on, and return the cursor to the origin.

We present two normal forms for an element $w \in L_2$, the {\em
right-first normal form} given by
$$rf(w)=a_{i_1} a_{i_2} \ldots a_{i_m} a_{-j_1} a_{-j_2} \ldots
a_{-j_l} t^{r}$$ and the {\em left-first normal form} given by
$$lf(w)=a_{-j_1} a_{-j_2} \ldots a_{-j_l} a_{i_1} a_{i_2} \ldots
a_{i_m} t^{r}$$ with $ i_m > \ldots i_2 > i_1 \geq 0 $ and $ j_l >
\ldots j_2 > j_1 > 0 $. In the successive prefixes of the right-first
normal form, the position of the cursor moves first toward the right,
and appropriate bulbs are illuminated in nonnegative positions.  In
later prefixes, the position of the cursor moves to the left, and once
its position is to the left of the origin, the appropriate bulbs in
negative positions are illuminated as well.  The left-first normal
form follows this procedure in reverse.

One or possibly both of these normal forms will lead to a minimal
length representative for $w \in L_2$ with respect to the wreath
product generating set, depending upon the final location of the
cursor relative to the origin. That position is easy to detect
from the sign of the exponent sum of $t$, given as $r$ above.  If
$r >0$, then the left-first normal form will produce a minimal length
representative for $w$, and if $r<0$, the right-first normal form
will produce a minimal length representative for $w$. If the
exponent sum of $t$ is zero, then both normal forms lead to
minimal length representatives for $w$.

Cleary and Taback \cite{\deadlamp} used these normal forms to compute
the word length of $w \in L_2$ with respect to the wreath product
generators as follows.
\begin{prop}[\cite{\deadlamp}, Proposition 3.2]
\label{D} Let $w \in L_2$ be in either normal form given
above, and define
$$D(w)=m+l+ \min\{2 j_l+i_m + | r-i_m|, 2 i_m+j_l+|r+j_l|\}.$$
The word length of $w$ with respect to the generating set
$\{a,t\}$ is exactly $D(w)$.
\end{prop}

The geodesic representatives for $w \in L_2$ arising from either
the left-first or right-first normal forms are not necessarily
unique.  Suppose that $\gamma$ is a geodesic word in $\{t,a\}$
representing $w \in L_2$, and $\gamma$ has two prefixes with the
same exponent sum.  This corresponds to a bulb at some position $k$
which is ``visited twice'' by the cursor during the construction of
$w$. If this bulb is illuminated in the word $w$, then there is a
choice as to whether it is illuminated during the first prefix
which leaves the cursor in this position, or the second.
Similarly, if the bulb in position $k$ is not illuminated in $w$,
then either it remains off during both prefixes, or is turned on
in the first, and off in the second.

If $w \in L_2$ has  illuminated bulbs which are visited more than once
 by the
cursor during the construction of $w$, then there will be multiple
possible geodesics representing $w$. A typical element $w$ where
the cursor is not left at the origin may have
$k$ illuminated bulbs which are visited exactly twice,  giving
$2^k$ possible
geodesics which represent $w$.  For elements where the final position of the
cursor is at the origin, there are  such families of geodesics from both the
right-first and left-first methods of construction. If bulb zero is illuminated and the cursor
remains at the origin with illuminated bulbs on both sides of the origin,
there will be 3 visits of the cursor to the origin, giving $3 \cdot 2^k$ possible
geodesic representatives for each of right-first and left-first manners of
constructing the element and thus $6 \cdot 2^k$ total possible 
geodesic representatives.
 It is not hard to see using
Proposition \ref{D} that all geodesic representatives for $w$
must be of this form.

\subsection{Normal forms for elements of $L_m$}

The normal forms given above for elements of $L_2$ have obvious
extensions to $L_m$ for $m>2$.  Occurrences of $a$ in the normal
forms above must now be replaced by $a^k$, for $k \in \{ -h, -h+1,
\cdots,-1,0,1,2, \cdots h \}$ where $h$ is the integer part of
$\frac{m}{2}$.  When $m$ is even, we may omit $a^{-h}$ to ensure
uniqueness, since $a^h=a^{-h}$ in $\Z_{2h}$.

There is an analogous definition for $D(w)$ in \cite{deadlamp}
when $w \in L_m$ which
again determines the word length of $w$ with respect to the wreath
product generating set. 

As with $L_2$, the left- and right-first normal forms for elements
of $L_m$ do not necessarily give all minimal length
representatives for group elements.  If a geodesic representative
for $w \in L_m$ has two prefixes with exponent sum $k$, and the
bulb in position $k$ in $w$ is illuminated to state $l<m$, then
there are many choices as to what state the bulb is left in at the
first prefix, with the second prefix allowing the bulb to be
illuminated to the final state.

\subsection{Geodesic languages for $L_m$}

We now prove that there is no collection of geodesic paths
representing elements of $L_m$ which is accepted by a finite state
automaton.  In this section we always consider $L_m$ with respect
to the wreath product generators $\{a,t\}$.

\begin{thm}\label{lem:LnotREG}
The lamplighter groups $L_m$ with respect to the wreath product
\gset  $\  \{a,t\}$ have no regular languages of geodesics.
\end{thm}

\textit{Proof}: 
Let $g_n=a_na_{-n} \in L_m$ be the group element corresponding to
a configuration of bulbs in which the two bulbs at distance $n$
from the origin are turned on to the first state, and the cursor
is at the origin.  This element has exactly two geodesic
representatives: $t^nat^{-2n}at^n$ and $t^{-n}at^{2n}at^{-n}$.

Suppose that there is a regular language of geodesics for $L_m$
with respect to this generating set.  Then the Pumping Lemma for
regular languages guarantees a pumping length $p$ for this
language.  Choose $n > p$ and consider $g_n$ as defined above.
Suppose that the first geodesic representative for $g_n$ is an
element of the regular language.

Since $n > p$, we can write $g_n=xyz$ such that $|xy|<p$, so
$x=t^i, y=t^j$ with $j>0$ and $j<p<n$, and
$z=t^{n-i-j}at^{-2n}at^n$. Then by the Pumping Lemma, $xy^2z$ must
also be in the language. So
$t^it^{2j}t^{n-i-j}at^{-2n}at^n=t^{n+j}at^{-2n}at^n$ would be
geodesic. This word has length $4n+j+2$ and corresponds to the
configuration of one lamp on at $n+j$, one lamp on at $-n+j$ and
the cursor at position $j$, to the right of the origin.  A shorter
word for this configuration is $t^{-n+j}at^{2n}at^{-n}$ which has
length $4n-j+2$, yielding a contradiction.

Similarly, if the second geodesic representative for $g_n$ was
part of the regular language, we would obtain an analogous
contradiction.  Thus neither geodesic representative for $g_n$ can
be part of any regular language of geodesics for $L_m$.
\hfill$\Box$

\begin{cor}\label{cor:lmconetypes}
The lamplighter group $L_m$ for $m \geq 2$ has infinitely many
cone types with respect to the wreath product \gset.
\end{cor}

\textit{Proof}: 
Theorem \ref{lem:LnotREG} states that the full language of
geodesics is not regular, so the contrapositive of Lemma
\ref{lem:coneimpliesreg} implies that the number of distinct cone
types is not finite.
\hfill$\Box$

One can observe this directly as well by considering the elements
$g_n=a_na_{-n}=t^nat^{-2n}at^n$ used in the proof of Theorem
\ref{lem:LnotREG}. The cone type of $t^nat^{-2n}a$ contains $t^n$
but not $t^{n+1}$, so for each $n$ we have a distinct cone type
with respect to $\{a,t\}$.

\begin{thm}\label{lem:lamp1counter}
There is a language of geodesics for $(L_m,\{a,t\})$ with a unique
representative for each element that is accepted by a 1-counter
automaton.
\end{thm}
\textit{Proof}: 
We first describe the 1-counter automaton accepting a language of
geodesics for $L_2$, and then give the generalization to $L_m$ for
$m > 2$. In $L_m$, for $m \geq 2$, each group element corresponds
to a configuration of bulbs in some states and a cursor position.
Such a configuration can always be obtained in the following
manner.

If the exponent sum of the generator $t$ in either normal form for
$w \in L_2$ is negative, so the cursor's final position is given by $i< 0$,
then we construct a minimal length representative for the element
in a ``right-first manner'' as follows.  Suppose that the rightmost
illuminated bulb in $w$ is in position $l$.  We begin the geodesic
representative with $t^{l}a$, which illuminates this bulb.  We
then add suffixes of the form $t^{-k}a$; each suffix puts the
cursor in the position of another bulb which must be turned on,
and this is accomplished via the generator $a$.  This is done
until the leftmost illuminated bulb is turned on.  We add a final
suffix of the form $t^n$ which brings the cursor to its position
in $w$.

\begin{figure}[ht!]
\psfrag{e}{$\epsilon$}
\psfrag{q0}{$q_0$}
\psfrag{q1}{$q_1$}
\psfrag{q2}{$q_2$}
\psfrag{A}{$A$}
\psfrag{e,+}{$(\epsilon,+)$}
\psfrag{e,-}{$(\epsilon,-)$}
\psfrag{a}{$a$}
\psfrag{t,+}{$(t,+)$}
\psfrag{T,-}{$(t^{-1},-)$}
\psfrag{ta,+}{$(ta,+)$}
\psfrag{Ta,-}{$(t^{-1}a,-)$}

  \bc
             \includegraphics[height=4cm]{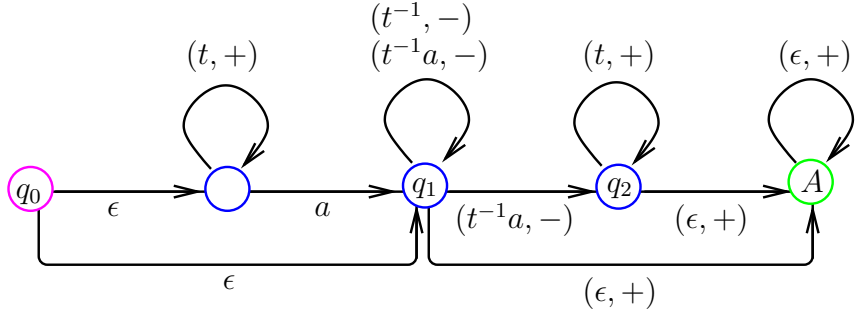}        
  \ec
  \caption{1-counter automaton accepting right-first normal form words for $L_2$.} 
  \label{fig:pic4A}
\end{figure}

The 1-counter automaton in Figure \ref{fig:pic4A}
 accepts this set of words. The counter keeps track of 
the current position of the cursor.

If the exponent sum of the generator $t$ in either normal form
given above is nonnegative, so that the cursor ends at a position
$i\geq 0$, then we construct a minimal length representative for
the element in a ``left-first manner''.  This is done by interchanging the
generators $t$ and $t^{-1}$ in the above method, and first turns
on the leftmost illuminated bulb.

\begin{figure}[ht!]
\psfrag{e}{$\epsilon$}
\psfrag{q0}{$q_0$}
 \psfrag{q3}{$q_3$}
\psfrag{q4}{$q_4$}
\psfrag{A}{$A$}
\psfrag{e,+}{$(\epsilon,+)$}
\psfrag{e,-}{$(\epsilon,-)$}
\psfrag{a}{$a$}
\psfrag{t,+}{$(t,+)$}
\psfrag{T,-}{$(t^{-1},-)$}
\psfrag{ta,+}{$(ta,+)$}
\psfrag{Ta,-}{$(t^{-1}a,-)$}

  \bc       
  \includegraphics[height=4cm]{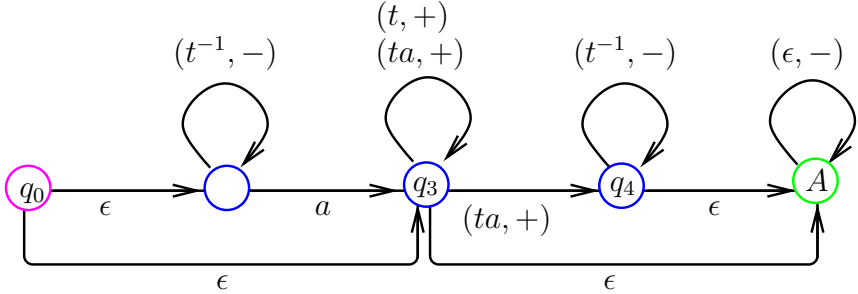}
  \ec
  \caption{1-counter automaton accepting left-first normal form words for $L_2$.} 
  \label{fig:pic4B}

\end{figure}

The 1-counter automaton in Figure  \ref{fig:pic4B} accepts
this set of words. Again the counter keeps track of 
the current position of the cursor.

This process gives a unique representative for each group element, 
and in \cite{\deadlamp} it is shown that these words are geodesic.

To construct a 1-counter automaton which accepts these
representatives of elements of $L_m$, we modify the automata
given in Figures \ref{fig:pic4A} 
and  \ref{fig:pic4B} as follows.  Additional
loops must be added to states $q_1$ and $q_3$, allowing for the bulbs
to be turned to any possible state.  Namely, there must be loops
at state $q_1$ with labels $(t^{-1}a^k,-)$ for  
$k \in \{ -h, -h+1,\cdots,-1,0,1,2, \cdots h \}$, where $h$ is the integer part of
$\frac{m}{2}$.  When $m$ is even, we omit $a^{-h}$ to ensure
uniqueness, since $a^h=a^{-h}$ in $\Z_{2h}$. We also need to add edges from $q_1$ to $q_2$ 
 with labels $(t^{-1}a^k,-)$ for  
$k \in \{ \pm 1,\pm 2, \cdots \pm h \}$.

The loops (and edges) added to state $q_3$ are of the form $(ta^k,+)$, with
the same restrictions on $k$.
\hfill$\Box$

Now we consider the full language of geodesics.   As described above,
there may be many possible geodesic representatives for an element.
Bulbs  may be visited twice or
even three times in the extreme case, when the cursor ends back at the
origin, with bulbs illuminated to the left and right. There may
be geodesic representatives which turn on such a bulb at any of those
opportunities, but no geodesic representative can change the state of
a particular bulb more than once in $L_2$.

Clearly, remembering which bulbs have already been turned on is not a
job for a finite state automaton, as proved above, but it is easy to
avoid switching on and off the same bulb by keeping track of
previous switchings using a stack.  This is the
key to the proof of the next theorem.

\begin{thm}[Context free full \lan\ for $(L_2,\{a,t\})$] \label{cflamp}
The full \lan\ of geodesics for the lamplighter group with the
wreath product \gset\ is \cf.
\end{thm}
\textit{Proof}: 
We  describe the complete set of geodesic words which represent
each group element $w \in L_2$, which we view as a configuration
of illuminated bulbs along with a position of the cursor.

If the cursor position in $w$ is $i\leq 0$ with no bulbs to the right
of position $0$ illuminated, then geodesic representatives for this element are of
the form $g_0g_{-1}\ldots g_{-m}(v)$ where $g_i$ is either $t^{-1} $
or $at^{-1}$ (which illuminates the bulb at position $i$), $v$ is
either empty or $ag'_{-m}g'_{1-m}\ldots g'_r$, where $-m\leq r\leq
0$ and

\begin{eqnarray*}
g'_i =
\begin{cases}
t & \mathrm{if} \;\; g_i=at^{-1} \\

at \;\; \mathrm{or} \;\; t & \mathrm{otherwise}.
\end{cases}
\end{eqnarray*}

\begin{figure}[ht!]
\psfrag{e}{$\epsilon$}
\psfrag{q0}{$q_0$}
\psfrag{A}{$A$}
\psfrag{a}{$a$}
\psfrag{e,push s}{$(\epsilon$,push $\$)$}
\psfrag{e,pop s}{$(\epsilon$,pop $\$)$}
\psfrag{T}{$t^{-1}$}
\psfrag{aT}{$at^{-1}$}
\psfrag{T,push 0}{$(t^{-1}$, push $0)$}
\psfrag{aT,push 1}{$(at^{-1}$, push $1)$}
\psfrag{t,pop 0}{$(t$, pop $0)$}
\psfrag{t,pop 1}{$(t$, pop $1)$}
\psfrag{ta,pop 0}{$(ta$, pop $0)$}

  \bc
               \includegraphics[height=4.66cm]{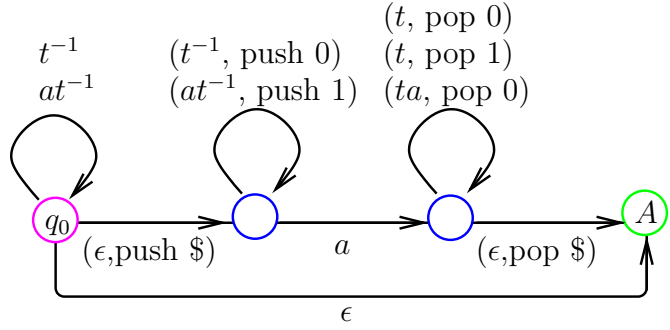} 

  \ec
  \caption{The pushdown automaton accepting the set of geodesics for
  elements of $L_2$ in which the cursor is in position $i\leq 0$ and
  there are no illuminated bulbs to the right of the origin.   The start state is
  labelled $q_0$.}
   \label{fig:pic5}
\end{figure}

We push a $1$ on the stack to indicate a bulb is switched on, and a $0$ if the bulb is not switched on in that position. Then when we return to a position we can only switch if a $0$ is on the top of the stack. The pushdown automaton in  Figure \ref{fig:pic5}
accepts these words.

If the cursor position in $w$ is $i\geq 0$ with no bulbs to the left
of position $0$ illuminated, then geodesic representatives for this element are of
the form $g_0g_1\ldots g_{n-1}(v)$ where $g_i$ is either $t $
or $at$ (if the lamp at position $i$ is illuminated), and $v$ is
either $\epsilon$ or $ag'_{n-1}\ldots g'_r$, where $0\leq r\leq n-1$ and

\begin{equation*}
g'_i =
\begin{cases}
t^{-1} & \mathrm{if} \;\; g_i=at \\

t^{-1}a \;\; \mathrm{or} \;\; t^{-1} & \mathrm{otherwise}.
\end{cases}
\end{equation*}

\begin{figure}[ht!]
\psfrag{e}{$\epsilon$}
\psfrag{q0}{$q_0$}
\psfrag{A}{$A$}
\psfrag{a}{$a$}
\psfrag{e,push s}{$(\epsilon$,push $\$)$}
\psfrag{e,pop s}{$(\epsilon$,pop $\$)$}
\psfrag{t}{$t$}
\psfrag{at}{$at$}
\psfrag{t,push 0}{$(t$, push $0)$}
\psfrag{at,push 1}{$(at$, push $1)$}
\psfrag{T,pop 0}{$(t^{-1}$, pop $0)$}
\psfrag{T,pop 1}{$(t^{-1}$, pop $1)$}
\psfrag{Ta,pop 0}{$(t^{-1}a$, pop $0)$}

  \bc

  \includegraphics[height=4.66cm]{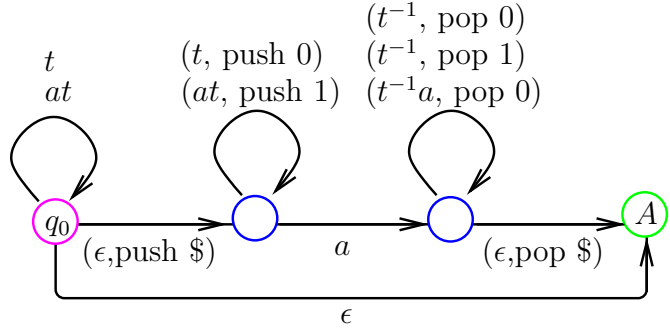} 

  \ec
  \caption{The pushdown 
  automaton accepting geodesics  in which the cursor is in position $i\geq 0$ and there are no
  illuminated bulbs to the left of the origin.  The start state is
  labelled $q_0$.}
   \label{fig:pic6}
\end{figure}

The pushdown automaton in the right of Figure \ref{fig:pic6}
accepts these words.

If the cursor position in $w$ is $i\leq 0$ with the rightmost
illuminated bulb at $n>0$, then geodesic representatives for this element are of the
form 
$$(u)f_1\ldots f_{n-1}tat^{-1}g_{n-1}\ldots g_1(v)g'_{-1}\ldots g'_{-m}(w)$$
 where $u$ is either $\epsilon$ or $a$ (turns bulb at $0$
on), $f_i$ is either $t$ or $ta$ (turns bulb at $i$ on),

\begin{equation*}
g_i =
\begin{cases}
t^{-1} & \mathrm{if} \;\; f_i=ta \\ 

at^{-1} \;\; \mathrm{or} \;\; t^{-1} & \mathrm{otherwise} \;\;(at^{-1}
\;\; \mathrm{turns}\;\; \mathrm{bulb}\;\; \mathrm{at}\;\; i \;\;
\mathrm{on}),
\end{cases}
\end{equation*}
\begin{equation*}
v =
\begin{cases}
a & \mathrm{if} \;\; u=\epsilon \;\; (\mathrm{turns}\;\; \mathrm{bulb}\;\;
\mathrm{at}\;\; i \;\; \mathrm{on}),\\

\epsilon & \mathrm{otherwise},
\end{cases}
\end{equation*}

$g'_i = t^{-1}$ or $t^{-1}a$ (turns the bulb at $i$ on), $w=\epsilon$
if the cursor ends at or to the left of the leftmost illuminated bulb,
or $t^{-1}atk_{-m}\ldots k_r$ with $r\leq 0$, for $-m\leq i<0$

 \begin{equation*}
k_i =
\begin{cases}
t & \mathrm{if} \;\; g'_i=t^{-1}a \\ 

ta \;\; \mathrm{or} \;\; t & \mathrm{otherwise} \;\;(at \;\;
\mathrm{turns}\;\; \mathrm{bulb}\;\; \mathrm{at}\;\; i \;\;
\mathrm{on}),
\end{cases}
\end{equation*}
and if r=0 then
\begin{equation*}
k_0 =
\begin{cases}
a & \mathrm{if} \;\; u=\epsilon \;\;\mathrm{and}\;\; v=\epsilon
\;\;(\mathrm{turns}\;\; \mathrm{bulb}\;\; \mathrm{at}\;\; i \;\;
\mathrm{on}),\\
\epsilon & \mathrm{otherwise}.
\end{cases}
\end{equation*}

\begin{figure}[ht!]
\psfrag{e}{$\epsilon$}
\psfrag{q0}{$q_0$}
\psfrag{A}{$A$}

\psfrag{e,push s}{$(\epsilon$, push $\$)$}
\psfrag{e,pop s}{$(\epsilon$, pop $\$)$}
\psfrag{e,pop h}{$(\epsilon$, pop $\#)$}
\psfrag{e,push h}{$(\epsilon$, push $\#)$}

\psfrag{a,push h}{$(a$, push $\#)$}
\psfrag{a,pop s}{$(a$, pop $\$)$}

\psfrag{taT}{$tat^{-1}$}
\psfrag{Tat}{$t^{-1}at$}

\psfrag{t,push 0}{$(t$, push $0)$}
\psfrag{ta,push 1}{$(ta$, push $1)$}
\psfrag{T,pop 0}{$(t^{-1}$, pop $0)$}
\psfrag{T,pop 1}{$(t^{-1}$, pop $1)$}
\psfrag{aT,pop 0}{$(at^{-1}$, pop $0)$}

\psfrag{T,push 0}{$(t^{-1}$, push $0)$}
\psfrag{Ta,push 1}{$(t^{-1}a$, push $1)$}
\psfrag{t,pop 0}{$(t$, pop $0)$}
\psfrag{t,pop 1}{$(t$, pop $1)$}
\psfrag{at,pop 0}{$(at$, pop $0)$}

  \begin{center}
               \includegraphics[height=4.66cm]{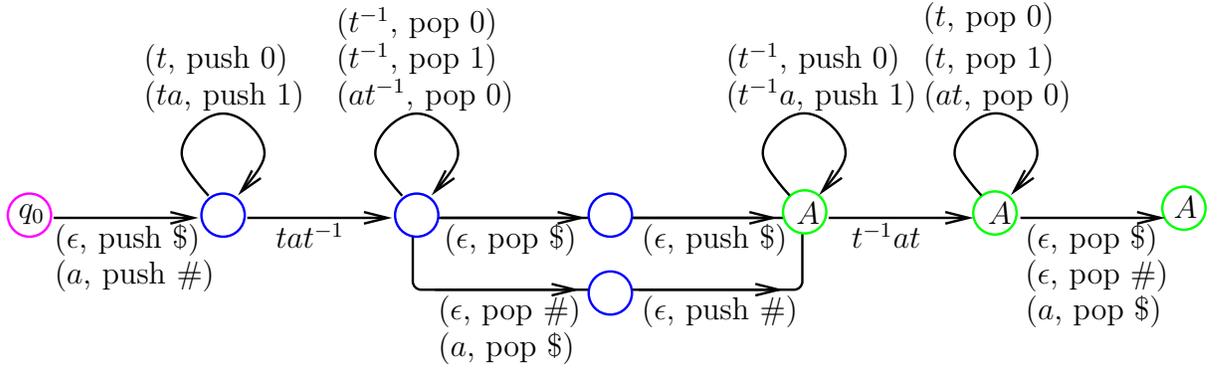}
  \end{center}
  \caption{The pushdown automaton accepting the set of geodesics for
  elements of $L_2$ in which the cursor is in position $i\leq 0$ and
  some bulb is illuminated at position $n>0$.  The state labelled $q_0$ is the
  start state.}
   \label{fig:pda-wreath1}
\end{figure}

As mentioned above, we can keep track of whether or not a bulb has
been illuminated as we move the cursor right then left then right using a
stack. In the pushdown automaton in Figure \ref{fig:pda-wreath1} we
push a $1$ on the stack to indicate a bulb is on, and a $0$ if the
bulb in that position is off.

We need to take care with the bulb at position 0, since there may be
three possible opportunities to illuminate it. We use different
``bottom of stack'' markers $\$$ and $\#$ depending on whether this
bulb has been turned on or not.

If the cursor position in $w$ is $i\geq 0$ with the leftmost
illuminated bulb at position $-m<0$, then geodesics for this element are of the
form \\
$(u)f_{-1}\ldots f_{1-m}t^{-1}atg_{1-m}\ldots g_{-1}(v)g'_1\ldots
g'_n(w)$ where $u$ is either $\epsilon$ or $a$ (turns bulb at $0$
on), $f_i$ is either $t^{-1}$ or $t^{-1}a$ (turns bulb at $i$ on),

 \begin{equation*}
g_i =
\begin{cases}
t & \mathrm{if} \;\; f_i=t^{-1}a \\ 

at \;\; \mathrm{or} \;\; t & \mathrm{otherwise},
\end{cases}
\end{equation*}
\begin{equation*}
v =
\begin{cases}
a & \mathrm{if} \;\; u=\epsilon \\

\epsilon & \mathrm{otherwise},
\end{cases}
\end{equation*}

$g'_i = t$ or $ta$ (turns the bulb at $i$ on), $w=\epsilon$
if the cursor ends at or to the right of the rightmost illuminated bulb,
or $tat^{-1}k_{n}\ldots k_r$ with $r\leq 0$, for $0<i\leq n$

\begin{equation*}
k_i =
\begin{cases}
t^{-1} & \mathrm{if} \;\; g'_i=ta \\ 

at^{-1} \;\; \mathrm{or} \;\; t^{-1} & \mathrm{otherwise},
\end{cases}
\end{equation*}
and if r=0 then
\begin{equation*}
k_0 =
\begin{cases}
a & \mathrm{if} \;\; u=\epsilon \;\;\mathrm{and}\;\; v=\epsilon \\

\epsilon & \mathrm{otherwise}.
\end{cases}
\end{equation*}

\begin{figure}[ht!]
\psfrag{e}{$\epsilon$}
\psfrag{q0}{$q_0$}
\psfrag{A}{$A$}

\psfrag{e,push s}{$(\epsilon$, push $\$)$}
\psfrag{e,pop s}{$(\epsilon$, pop $\$)$}
\psfrag{e,pop h}{$(\epsilon$, pop $\#)$}
\psfrag{e,push h}{$(\epsilon$, push $\#)$}

\psfrag{a,push h}{$(a$, push $\#)$}
\psfrag{a,pop s}{$(a$, pop $\$)$}

\psfrag{taT}{$tat^{-1}$}
\psfrag{Tat}{$t^{-1}at$}

\psfrag{t,push 0}{$(t$, push $0)$}
\psfrag{ta,push 1}{$(ta$, push $1)$}
\psfrag{T,pop 0}{$(t^{-1}$, pop $0)$}
\psfrag{T,pop 1}{$(t^{-1}$, pop $1)$}
\psfrag{aT,pop 0}{$(at^{-1}$, pop $0)$}

\psfrag{T,push 0}{$(t^{-1}$, push $0)$}
\psfrag{Ta,push 1}{$(t^{-1}a$, push $1)$}
\psfrag{t,pop 0}{$(t$, pop $0)$}
\psfrag{t,pop 1}{$(t$, pop $1)$}
\psfrag{at,pop 0}{$(at$, pop $0)$}

  \begin{center}
               \includegraphics[height=4.66cm]{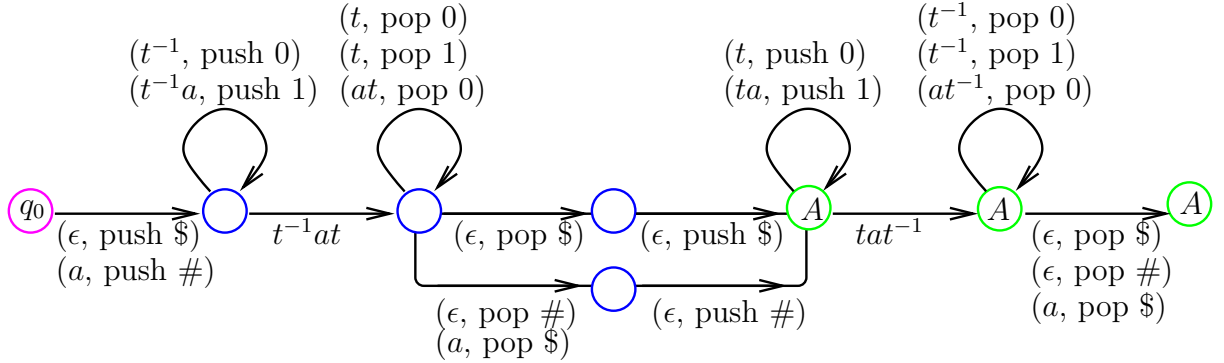}
  \end{center}
  \caption{The pushdown automaton accepting the set of geodesics for
  elements of $L_2$ in which the cursor is in position $i\geq 0$ and some bulb is illuminated at position $n<0$.  The
  state labelled $q_0$ is the start state.}
   \label{fig:pda-wreath2}
\end{figure}

It follows from \cite{\deadlamp} that all possible geodesics for
elements of $L_2$ in this generating set have one of these forms, and
is accepted by one of the pushdown automata in Figures 
\ref{fig:pic5} - \ref{fig:pda-wreath2}.  By Lemma
\ref{lem:CFunion} the union of these four languages is \cf.
\hfill$\Box$

Next, we show that the full language of geodesics is not counter.  The
proof mimics Elder's proofs in \cite{\MeBSCF} which show that there are
\cf\ languages that are not counter, and is derived from the fact that
one can write out a string on three letters of arbitrary length that
has no concurrent repeating subwords, due to Thue and Morse
\cite{\MeBSCF}. We call such a string of letters with no repeating
subwords a {\em Thue-Morse word}.

\begin{thm}
\label{thm:notcounter} The \lan\ of all geodesics for the
lamplighter group $L_m$ with the wreath product \gset\ $\{a,t\}$
is not counter.
\end{thm}
\textit{Proof}: 
Suppose that $C$ is the full language of geodesics for $L_m$, and
that it is counter. As earlier, we let $h$ be the integer part of
$\frac{m}{2}$, and consider the regular language $L$ of all
strings consisting of $a^h$ and $t^{\pm 1}$, beginning with $t$
and ending with $a$, with at most $3$ consecutive $t^{\pm 1}$
symbols. The intersection of $C$ with $L$ is counter by Lemma
\ref{lem:closurekcounter}.

We form a language $E$ on the letters $\pm 1,\pm 2 , \pm 3$ to encode
words from $C\cap L$ as follows.  Encode a word
$t^{w_1}a^ht^{w_2}a^h\ldots t^{w_n}a$ by its $t$-exponents, to obtain
$w_1w_2\ldots w_n$.  So for example the word
$w=t^2a^ht^3a^hta^ht^2a^h$ is encoded as $e=2312$.

An encoded word $e=e_1e_2 \cdots e_k$ with all $e_i$ positive
represents a configuration of bulbs with some subset of the bulbs
in positions from $1$ and $\sum e_i$ illuminated to state $h$. We
consider a Thue-Morse word $w$ in $L \cap C$, which has by
definition no repeated subwords, and encode it to form a string in
$E$ which also has no repeated subwords.  This word and its
encoding both represent the same configuration of bulbs, where all
of the illuminated bulbs are in positions between $1$ and $\sum
e_i$.

We append a suffix to $w$ as follows, so that the new word has all
bulbs in positions between 1 and $\sum e_i$ illuminated to state
$h$. Namely, we attach strings of the form $t^{-s}a^h$ to the end
of $w$, for $s \in \{1,2,3\}$, until all bulbs are illuminated to
state $h$.

For example, if $m=2$ and $e = 2312$ represents a configuration of
illuminated bulbs in positions from $1$ to $8$, then
$2312(-1)(-3)(-1)(-2)$ is the encoding for the word with all bulbs
illuminated between these positions.

When $e$ is a Thue-Morse word, there is a nice rewriting of the
positive part of the word to obtain the negative part, but for
this argument we merely need to observe the following.  If $e$ is
a Thue-Morse word, then we let $f$ be the encoding suffix described above
corresponding to the group element $v$, so that $wv$ represents a
group element with all bulbs illuminated to state $h$ from
positions 1 to $\sum e_i$ with the cursor at the origin. Suppose
that $w'$ represents a different configuration of bulbs
illuminated to state $h$ at positions between $1$ and $\sum e_i$.
Then $w'v$ represents a word in which some bulb is turned to state
$2h$.  If $2h = m$, then this bulb has been turned on and then
off, thus the word is not geodesic. If $2h = m-1$, then there is a
bulb which is turned to state $2h$; a geodesic representative for
this element would use a single generator, $a^{-1}$, to achieve
this state. The point is that the pairs of positive and negative
words defined here are carefully chosen to interleave illumination
of the bulbs, thus ensuring that changing
either of the parts alone will render the word non-geodesic.

If the language $C \cap L$ is counter, then it is accepted by some
counter automaton.  We construct a counter automaton accepting the
encoded language $E$ as follows. For every path in the finite
state automaton labelled by $t^ia^h$ plus some counters $u$, where
$i\in \{\pm1,\pm2,\pm3\}$, replace the path by an edge labelled
$(i,u)$. This gives a new counter automaton with the same
counters, accept states and start state, accepting the language of
encoded words. Let $p$ be the swapping length for this language,
guaranteed by the Swapping Lemma (Lemma \ref{lem:swap}).

Take a Thue-Morse word in $1,2,3$ of length greater than $2p+1$,
then append a word in $-1,-2,-3$ so that the full word represents
a group element with all lamps illuminated in positions from $1$
to some large positive integer.  Call this element $w$.

By the Swapping Lemma, we must be able to swap two adjacent
subwords in the initial positive segment of $w$ and obtain another
word $w'$ in the language.  But if we do this the negative part of
$w'$ is the same as the negative part of $w$ and is now no longer
compatible with the changed positive part. Thus we have an
encoding of a word that is not geodesic, so $w'$ cannot be in $C
\cap L$, a contradiction.  Thus, the full language of geodesics is
not counter.
\hfill$\Box$

\section{Geodesic languages for the lamplighter groups with the automata \gset}
\label{sec:lampautomata}

Grigorchuk and Zuk \cite{\Grig} prove that the lamplighter group $L_2$
is an example of an automata group.  The natural generating set which
arises from this interpretation is $\{t, ta\}$ which we will call the
automata generating set for $L_2$.  They compute the spectral radius
of $L_2$ with respect to this generating set and find remarkably that
it is a discrete measure.  Bartholdi and Sunik \cite{bs} prove that all lamplighter groups $L_m$ are automata groups.

We now compare the language theoretic properties of $L_2$ with respect to the automata generating set $\{t, ta\}$ to the analogous properties described in the previous sections.  We still view elements of $L_2$ as a configuration of light bulbs,
with a cursor pointing to an integral position.  However, the
generator $ta$ combines the two basic ``motions'' of the wreath
product generators, namely multiplication by $ta$ both moves the
cursor and turns on a bulb. We show here that, as with respect to
the wreath product generating set, there is no regular language of
geodesics for $L_2$ and there is a counter language of unique
geodesic representatives.  In contrast to the  wreath product
generating set, however,  the full language of geodesics with respect to the
automata generating set is 1-counter.  

\subsection{Geodesic paths with respect to the automata generators}

We can construct minimal length representatives with respect to
the automata generating set using the same normal forms for
elements of $L_2$ described above when considering the wreath
product generating set, as described by Cleary and Taback \cite{\lamptat}.

We state an analogue of Proposition \ref{D} which will allow us to
recognize which of these normal forms are geodesic with respect to
the automata generating set.  As described in \cite{\lamptat}, the
length of an element with respect to the automata generating set
depends only upon the positions of the leftmost- and
rightmost-illuminated bulbs and the final location of the cursor,
since it is possible to turn on or off intermediate bulbs with no
additional usage of generators by choosing to use $ta$ instead of
$t$ for moving the cursor.

\begin{defn}
Let $w = a_{i_1} a_{i_2} \ldots a_{i_m} a_{-j_1} a_{-j_2} \ldots
a_{-j_l} t^{r} \in L_2$, with $0 < i_1 < i_2 \cdots < i_m $ and $0
\leq j_1 < j_2 \cdots < j_l$. If $l=0$, there are no bulbs
illuminated at or to the left of the origin and we set $D'(w)=
i_m+ | r -i_m|$. Otherwise, we set
 $$D'(w)= \min\{ 2 (j_l+1)+
 i_m+|r-i_m|,2 i_m+ j_l + 1 + |r+j_l+1|\}.$$
\end{defn}

With respect to the automata generating set it is more convenient
to group the bulb in position $0$ with the bulbs in negative
positions in the normal form given above than with the bulbs in
positive positions, as in \cite{\lamptat}.  This is done because
the element $a \in L_2$ has length two with respect to this
generating set, and is explained fully in \cite{\lamptat}.

\begin{prop}[\cite{\lamptat}, Proposition 2.4]
\label{prop:D} The word length of $w \in L_2$ with respect to
 the automata generating set
$\{t,ta\}$ is given by $D'(w)$.
\end{prop}

We use Proposition \ref{prop:D} to show that $L_2$ has no regular
language of geodesics with respect to the automata generating set.

\begin{thm}\label{lem:LnotREGauto}
The lamplighter group with the automata \gset\ has no regular
language of geodesics.
\end{thm}

\textit{Proof}:  The argument is similar to that for the wreath product
generating set.  Again, we consider words of the form
$g_n=a_na_{-n}$. This element has length $4n+2$ with respect to
the automata generating set.  There are again two families of
geodesic representatives for $g_n$, those arising from the
right-first normal forms, such as $t^n (ta)^{-1} t^{-2n} (ta)
t^n$, and those arising from the left-first normal forms, such as
$t^{-n-1} (ta) t^{2n} (ta)^{-1} t^{-n+1}$.  Applying the Pumping Lemma as in the proof of Theorem \ref{lem:LnotREG} yields the desired contradiction.

%
%
%
\hfill$\Box$

As before, it follows from Lemma \ref{lem:coneimpliesreg} that
there are infinitely many cone types.

\begin{cor}\label{cor:l2conetypes}
The lamplighter group has infinitely many cone types with respect
to the automata \gset.
\end{cor}

Again, we can observe this directly by considering  the
elements\\ $g_n=a_na_{-n}=t^n (ta)^{-1} t^{-2n-1} (ta)
t^n$. The cone type of $g_n t^{-k} $ contains $t^k$ but not
$t^{k+1}$ for $0 \leq k \leq n$, so for each $k$ we have a
distinct cone type and as $n$ increases,  we have infinitely
many cone types.

Though the languages of geodesics with respect to the automata
generating set cannot be regular, there are geodesic languages with respect
to this generating set which are 1-counter.  
In contrast to the the wreath
product generating set, it is now possible to have geodesic words which
change the state of a bulb several times, so it is not essential to
remember which bulbs have already been illuminated by using a stack.

\begin{thm}\label{lem:lamp1counterauto}
There is a language of geodesics for $(L_2,\{t,(ta)\})$ with a
unique representative for each element that is accepted by a
1-counter automaton.
\end{thm}
\textit{Proof}: 
Let $w \in L_2$.  We choose the unique representative for $w$ in our
language to be of one of the following forms.

Suppose that $w$ has illuminated bulbs to the right of the origin and
also at or to the left of the origin and the position of the cursor is
at or to the left of the origin.  Let the position of the rightmost
illuminated bulb be $n >0$, and the position of the leftmost
illuminated bulb be $-m$ with $m \leq 0$.  We choose a representative
of $w$ of the form $t^n (ta)^{-1} g_{n-1} g_{n-2} \ldots g_{1-m}
(ta)^{-1} t^l$ where $l \leq m$ and

\begin{equation*}
g_i =
\begin{cases}
(ta)^{-1} & \mathrm{if} \;\; \mathrm{bulb}\;\; i \;\; \mathrm{is} \;\;
\mathrm{on}\\ 

t^{-1} & \mathrm{if} \;\; \mathrm{bulb}\;\; i \;\;\mathrm{is} \;\;
\mathrm{off}.
\end{cases}
\end{equation*}

\begin{figure}
\psfrag{e}{$\epsilon$}
\psfrag{q0}{$q_0$}
\psfrag{A}{$A$}
\psfrag{e,+}{$(\epsilon,+)$}
\psfrag{t}{$t$}
\psfrag{t,+}{$(t,+)$}
\psfrag{T,-}{$(t^{-1},-)$}
\psfrag{(TA)}{$(ta)^{-1}$}
\psfrag{(TA),-}{$((ta)^{-1},-)$}

 \bc
              \includegraphics[height=4cm]{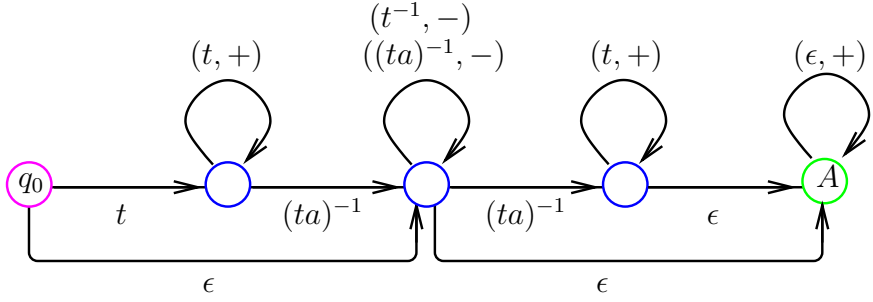}
  \ec
\caption{A $1$-counter automaton accepting geodesics representing elements $w$ with illuminated bulbs at $n>0$ and $-m$ with $m\leq 0$
and cursor at $k\leq 0$.}\label{fig:1counter-tta-gset}
\end{figure}

The 1-counter automaton must check that the number of $t$ letters at
the start and end of a word is at most the number of
$(ta)^{-1},t^{-1}$ letters, and the word starts with at least one $t$
letter. See Figure \ref{fig:1counter-tta-gset}.

If there are illuminated bulbs only at or to the left of the origin,
then we choose a representative for $w$ of the form $g_0 g_{-1} \ldots
g_{-m}(v)$ where $g_i$ is either $t^{-1}$ or $(ta)^{-1}$ as defined
above, and $v$ is either empty (if the cursor is to the left of the
leftmost illuminated bulb), or $(ta)^{-1}t^l$ where $1\leq l\leq
m$. These words are obtained in the 1-counter automaton in Figure
\ref{fig:1counter-tta-gset}, by following the $\epsilon$ edge from the
start state.

Suppose that $w$ has illuminated bulbs to the right of the origin and
also at or to the left of the origin and the position of the cursor is
to the right of the origin.  Let $n$ and $m$ be as above, and choose a
representative of the form $t^{-m-1} (ta) g_{1-m} \ldots g_{n-1} (ta)
t^l$ where and $l> - n$ and

\begin{equation*}
g_i =
\begin{cases}
(ta) & \mathrm{if} \;\; \mathrm{bulb}\;\; i \;\; \mathrm{is} \;\;
\mathrm{on}\\ 

t & \mathrm{if} \;\; \mathrm{bulb}\;\; i \;\;\mathrm{is} \;\;
\mathrm{off}.

\end{cases}
\end{equation*}

\begin{figure}
\psfrag{e}{$\epsilon$}
\psfrag{q0}{$q_0$}
\psfrag{A}{$A$}
\psfrag{e,-}{$(\epsilon,-)$}
\psfrag{t}{$t$}
\psfrag{t,+}{$(t,+)$}
\psfrag{T}{$t^{-1}$}
\psfrag{T,-}{$(t^{-1},-)$}
\psfrag{(ta)}{$(ta)$}
\psfrag{(ta),+}{$((ta),+)$}

 \bc
              \includegraphics[height=4cm]{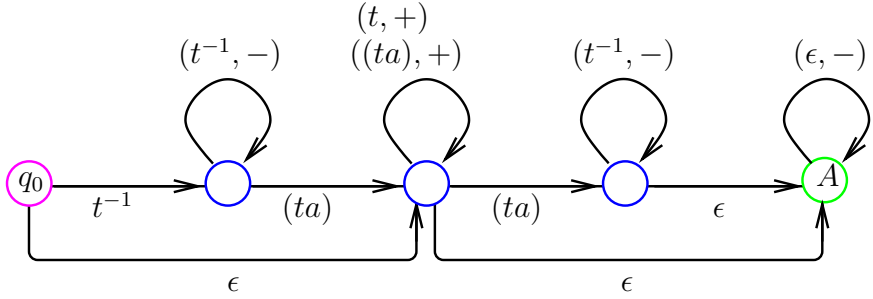}
  \ec
\caption{A $1$-counter automaton accepting geodesics representing elements $w$ with illuminated bulbs at $n>0$ and $-m$ with $m\leq 0$
and cursor at $k> 0$.}

\label{fig:1counter-tta-gsetPOS}.
\end{figure}

The 1-counter automaton for these words must check that the number of
$t^{-1}$ letters at the start and end of a word is at most the total number of
$(ta),t$ letters, and the word starts with at least one $t^{-1}$
letter. See Figure \ref{fig:1counter-tta-gsetPOS}.

If there are illuminated bulbs only to the right of the origin, then
we choose a representative for $w$ of the form $(u)g_1 g_{2} \ldots
g_{n-1}(v)$ where $u$ is either empty (if the bulb at the origin is
off) or $t^{-1}(ta)$ (if the bulb at the origin if on), $g_i$ is
either $t$ or $(ta)$ as defined above, and $v$ is either empty (if the
cursor is at or to the right of the rightmost illuminated bulb), or
$(ta)t^{-l}$ where $1\leq l\leq n$. These words are obtained in the
1-counter automaton in Figure \ref{fig:1counter-tta-gsetPOS}, by
following the $\epsilon$ edge from the start state.
\hfill$\Box$

Again, in contrast to the wreath product generating set, the
full language of geodesics with respect to the automata generating
set can be recognized by a simpler machine-- it does not require a stack to keep
track of which bulbs have been  illuminated.

\begin{thm}[1-counter full \lan\ for $(L_2,\{t,(ta)\})$]
\label{thm:tta1countfull}
The full \lan\ of geodesics for the lamplighter group with the
automata \gset\ is accepted by a 1-counter automaton.
\end{thm}
\textit{Proof}: 
We exhibit four 1-counter automata, each recognizing a type of
geodesic word.  There are several variations within each family of
geodesics which must be recognized as well.

We first consider geodesic representatives for elements $w \in L_2$ in
which all illuminated bulbs lie to the right of the origin, with the
rightmost illuminated bulb in position $n > 0$.

In this case, we build an automaton which recognizes geodesics of
the following forms.

\bi
\item If the cursor position in $w$ is at or to the right of $n$, then
geodesics are words of the form $g_1 \ldots g_n$ with $g_i=t$ or
$(ta)$.
\item If the cursor is to the left of $n$, then geodesics are words of
the form $g_1\ldots g_{n-1}uv$ with $g_i=t$ or $(ta)$, $u$ is
$(ta)t^{-1}$ or $t(ta)^{-1}$, $v$ is $\epsilon$, $g'_{n-1}\ldots g'_r$
with $r\geq 0$, or $g'_{n-1}\ldots g'_0t^{-r}$ with $r>0$, and
$g'_i=t^{-1}$ or $(ta)^{-1}$.  \ei

\begin{figure}
\psfrag{e}{$\epsilon$}
\psfrag{q0}{$q_0$}
\psfrag{A}{$A$}
\psfrag{e,-}{$(\epsilon,-)$}
\psfrag{t}{$t$}
\psfrag{T}{$t^{-1}$}
\psfrag{t,+}{$(t,+)$}
\psfrag{T,-}{$(t^{-1},-)$}
\psfrag{(TA),-}{$((ta)^{-1},-)$}
\psfrag{t(TA),-}{$t(ta)^{-1}$}
\psfrag{(ta)T}{$(ta)t^{-1}$}
\psfrag{t(TA)}{$t(ta)^{-1}$}
\psfrag{(ta)}{$(ta)$}
\psfrag{(ta),+}{$((ta),+)$}

 \bc
  \includegraphics[height=4cm]{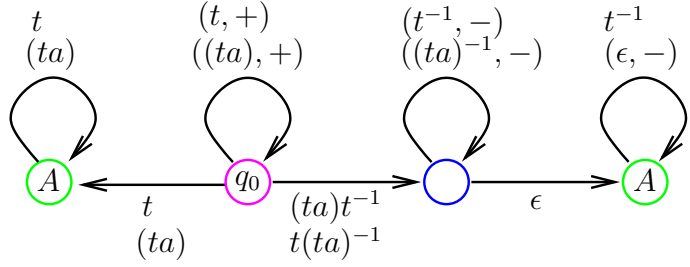}
   \ec
         \caption{A $1$-counter automaton accepting geodesics representing elements $w$ with all illuminated bulbs to the right of the origin.}
\label{fig:full-1countA}
\end{figure}

These forms are the language of the counter automaton in Figure
\ref{fig:full-1countA}. The first type are obtained by following the
$t$ or $(ta)$ edge to the  accept state on the left, and no counters are needed.

Second, we consider geodesic representatives for elements $w$ in
which all illuminated bulbs lie at or to the left of the origin, with
the leftmost illuminated bulb in position $-m$ with $m \geq 0$.

In this case, we build an automaton which recognizes geodesics of
the following forms.

\bi
\item If the cursor position in $w$ is to the left of $-m$, then
geodesics are words of the form $g_0 \ldots g_{-m}$ with $g_i=t^{-1}$
or $(ta)^{-1}$.
\item If the cursor is at or to the right of $-m$, then geodesics are
words of the form $g_0\ldots g_{1-m}uv$ with $g_i=t^{-1}$ or
$(ta)^{-1}$, $u$ is $(ta)^{-1}t$ or $t^{-1}(ta)$, $v$ is $\epsilon$,
$g'_{1-m}\ldots g'_r$ with $r\leq 0$, or $g'_{1-m}\ldots g'_0t^r$ with
$r>0$, and $g'_i=t$ or $(ta)$.  \ei

These forms are the language of the counter automaton in Figure
\ref{fig:full-1countA}. The first type are obtained by following the
$t^{-1}$ or $(ta)^{-1}$ edge to the accept state on the left, and no counters are
needed.

\begin{figure}
\psfrag{e}{$\epsilon$}
\psfrag{q0}{$q_0$}
\psfrag{A}{$A$}
\psfrag{e,+}{$(\epsilon,+)$}
\psfrag{T}{$t^{-1}$}
\psfrag{t}{$t$}
\psfrag{t,+}{$(t,+)$}
\psfrag{T,-}{$(t^{-1},-)$}
\psfrag{(TA),-}{$((ta)^{-1},-)$}
\psfrag{(TA)}{$(ta)^{-1}$}
\psfrag{(TA)t}{$(ta)^{-1}t$}
\psfrag{T(ta)}{$t^{-1}(ta)$}
\psfrag{(ta),+}{$((ta),+)$}

 \bc
  \includegraphics[height=3.33cm]{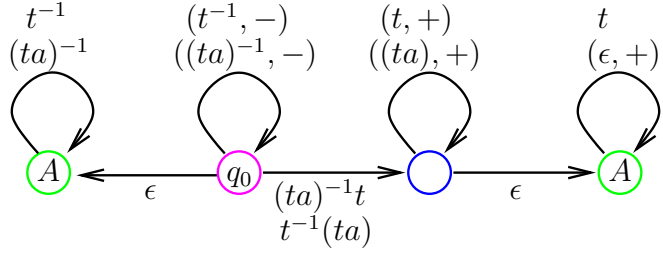} 
  \ec
         \caption{A $1$-counter automaton accepting geodesics representing elements $w$ with all illuminated bulbs at or to the left of the origin.}
\label{fig:full-1countB}
\end{figure}

Finally, we consider geodesic representatives for elements $w$ in
which illuminated bulbs lie to the right and at or to the left of the
origin, with the rightmost illuminated bulb at $n>0$, and the leftmost
illuminated bulb at $-m$, for $m \leq 0$.

In this case, we build an automaton which recognizes geodesics of
the following forms.

\bi
\item If the cursor position in $w$ is at or to the left of the origin, then
geodesics are words of the form 
 $g_1\ldots g_{n-1}ug'_{n-1}\ldots g'_0  v $ with $g_i=t$ or $(ta)$, $u$ is
$(ta)t^{-1}$ or $t(ta)^{-1}$, $g'_i=t^{-1}$ or $(ta)^{-1}$,   $v$ is
$\epsilon$ or $h_1\ldots h_rxy$ with 
  with $h_i=t^{-1}$ or $(ta)^{-1}$, $x$ is
$(ta)^{-1}t$ or $t^{-1}(ta)$, and $y$ is $h'_r\ldots h'_{m}$
 with $h'_i=t$ or $(ta)$ and $m\geq 0$.

\item If the cursor is at or to the right of the origin, then geodesics are
words of the form 
 $g_1\ldots g_{n-1}ug'_{n-1}\ldots g'_0  v $ with $g_i=t^{-1}$ or $(ta)^{-1}$, $u$ is
$(ta)^{-1}t$ or $t^{-1}(ta)$, $g'_i=t$ or $(ta)$,   $v$ is
$\epsilon$ or $h_1\ldots h_rxy$ with 
  with $h_i=t$ or $(ta)$, $x$ is
$(ta)t^{-1}$ or $t(ta)^{-1}$, and $y$ is $h'_r\ldots h'_{m}$ with $h'_i=t^{-1}$ or $(ta)^{-1}$
 and $m\geq 0$.

\ei

\begin{figure}

\psfrag{e}{$\epsilon$}
\psfrag{q0}{$q_0$}
\psfrag{A}{$A$}
\psfrag{e,+}{$(\epsilon,+)$}
\psfrag{T}{$t^{-1}$}
\psfrag{t}{$t$}
\psfrag{t,+}{$(t,+)$}
\psfrag{T,-}{$(t^{-1},-)$}
\psfrag{(TA),-}{$((ta)^{-1},-)$}
\psfrag{(TA)}{$(ta)^{-1}$}
\psfrag{(TA)t}{$(ta)^{-1}t$}
\psfrag{T(ta)}{$t^{-1}(ta)$}
\psfrag{t(TA)}{$t(ta)^{-1}$}
\psfrag{(ta)T}{$(ta)t^{-1}$}
\psfrag{(ta),+}{$((ta),+)$}

 \bc
   \includegraphics[height=4.33cm]{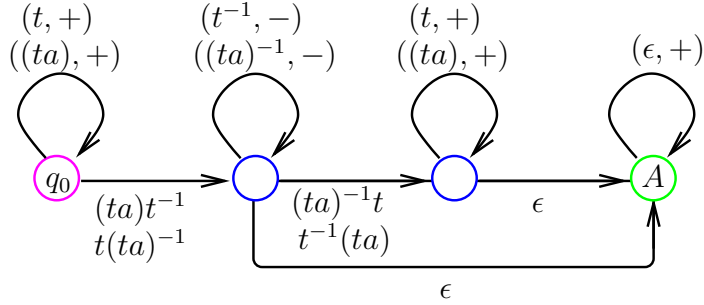} 

 \ec
         \caption{1-counter automaton for geodesics with bulbs to the
         right and at or to the left of the origin and cursor at or to the left of the origin.}
\label{fig:pic13}
\end{figure}

\begin{figure}

\psfrag{e}{$\epsilon$}
\psfrag{q0}{$q_0$}
\psfrag{A}{$A$}
\psfrag{e,-}{$(\epsilon,-)$}
\psfrag{T}{$t^{-1}$}
\psfrag{t}{$t$}
\psfrag{t,+}{$(t,+)$}
\psfrag{T,-}{$(t^{-1},-)$}
\psfrag{(TA),-}{$((ta)^{-1},-)$}
\psfrag{(TA)}{$(ta)^{-1}$}
\psfrag{(TA)t}{$(ta)^{-1}t$}
\psfrag{T(ta)}{$t^{-1}(ta)$}
\psfrag{t(TA)}{$t(ta)^{-1}$}
\psfrag{(ta)T}{$(ta)t^{-1}$}
\psfrag{(ta),+}{$((ta),+)$}

 \bc
   \includegraphics[height=4.33cm]{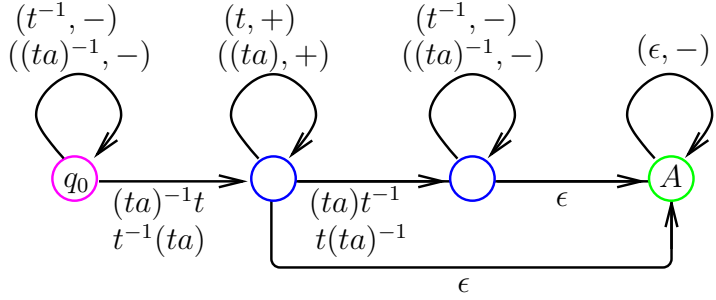} 

 \ec
         \caption{1-counter automaton for geodesics with bulbs to the
         right and at or to the left of the origin and cursor at or to the right of the origin.}
\label{fig:pic14}
\end{figure}

Each of these forms are accepted by the counter automata in Figures
\ref{fig:pic13} and \ref{fig:pic14}.\hfill$\Box$

Theorems \ref{thm:notcounter} and \ref{thm:tta1countfull} together
show that $L_2$ has at least one generating set which yields a full
language of geodesics that is counter, and one generating set whose
full language of geodesics is not counter.

\section{Geodesic languages for Thompson's group $F$}
\label{sec:F}

Thompson's group $F$ is a fascinating group which is studied from
many different perspectives.
Analytically, $F$ is understood as a group of piecewise linear
homeomorphisms of the interval, whose finitely-many discontinuities of slope have dyadic
rational coordinates, and whose linear pieces have slopes which
are powers of $2$. See Cannon, Floyd and Parry \cite{\cfp} for an
excellent introduction to $F$.

Algebraically, $F$ has two standard presentations: an infinite one
with a convenient set of normal forms, and a finite one, with respect
to which we study the metric. These presentations are as follows: 
\begin{eqnarray*}
F & = & \langle x_0,x_1
|[x_0x_1^{-1},x_0^{-1}x_1x_0],[x_0x_1^{-1},x_0^{-2}x_1x_0^2] \rangle \\
& = & 
 \langle x_i, i \geq 0 | x_i^{-1} x_j x_i = x_{j+1} , \textrm{for } i < j
\rangle.
\end{eqnarray*}

Geometrically, $F$ is studied as a group of pairs of finite binary
rooted trees, where group multiplication is analogous to function
composition.   Fordham \cite{\blakegd} developed a remarkable
way of measuring the word length of an element of $F$, with
respect to the finite generating set $\{x_0,x_1\}$, directly from
the tree pair diagram representing the element.  Belk and Brown
\cite{\belkbrown} and Guba \cite{\gubagrowth}
 also have geometric methods
for computing the word length of an element of $F$ in this
generating set.

All of these methods for computing word length in the generating
set $\{x_0,x_1\}$ have led to greater understanding of the the
geometry of the Cayley graph of $F$ with respect to this
generating set. For example, Cleary and Taback \cite{\ct} show that
this Cayley graph is not almost convex. Belk and Bux
\cite{\belkbux} show that this Cayley graph is additionally not
minimally almost convex. Burillo \cite{\josegrowth}, Guba
\cite{\gubagrowth}, and Belk and Brown \cite{\belkbrown} have studied
the growth of the group by trying to compute the size of balls in
this Cayley graph.  These volumes have been estimated, but  the
exact growth function of $F$ is not known and it is not even known
if it is rational.

In this section, we prove directly that $F$ contains infinitely
many cone types with respect to the standard finite generating set
$\{x_0,x_1\}$. This fact also follows as a corollary to Corollary
\ref{cor:conetypes} below.  Our explicit example uses a family of
{\em seesaw elements}, described by Cleary and Taback \cite{seesaw},
which are words for which there are two
different possible suffixes for geodesic representatives.

We extend the results of Section \ref{sec:F} to general groups containing seesaw elements in Section \ref{sec:generalcase}.

\begin{defn}
An element $w$ in a finitely generated group $G$ with finite
generating set $X$ is a {\em seesaw element of swing $k$ with respect
to a generator $g$} if the following conditions hold.
\bi
\item  Right multiplication by both $g$ and $g^{-1}$ reduces the
word length of $w$; that is, $|wg^{\pm 1}| = |w| - 1$, and for all
$h \in X  \smallsetminus \left\{g^{\pm1}\right\}$, we have
$|wh^{\pm 1}| \geq |w|$.

\item Additionally, $|wg^l| = |wg^{l-1}| - 1$   for integral $l
\in [1,k]$, and $|wg^{m}h^{\pm 1}| \geq |wg^{m}|$ for all $h
\in X  \smallsetminus \left\{g\right\}$ and integral $m \in
[1,k-1]$.

\item Similarly,  $|wg^{-l}| = |wg^{-l+1}| - 1$  for integral $l
\in [1,k]$, and $|wg^{-m}h^{\pm 1}| \geq |wg^{-m}|$ for all $h
\in X  \smallsetminus \left\{g^{-1}\right\}$  for integral $m \in
[1, k-1]$.
\ei
\end{defn}

Seesaw elements pose difficulty for finding unique geodesic
representatives consistently for group elements, and are used to show
that $F$ is not combable by geodesics in \cite{\seesaw}.  We notice
from the definition that if $w$ is a seesaw element of swing $k$ then
geodesic representatives for $w$ can have exactly two suffixes of
length $k$-- either $g^k$ or $g^{-k}$.

\begin{figure}
\bc

\includegraphics[width=14.2cm]{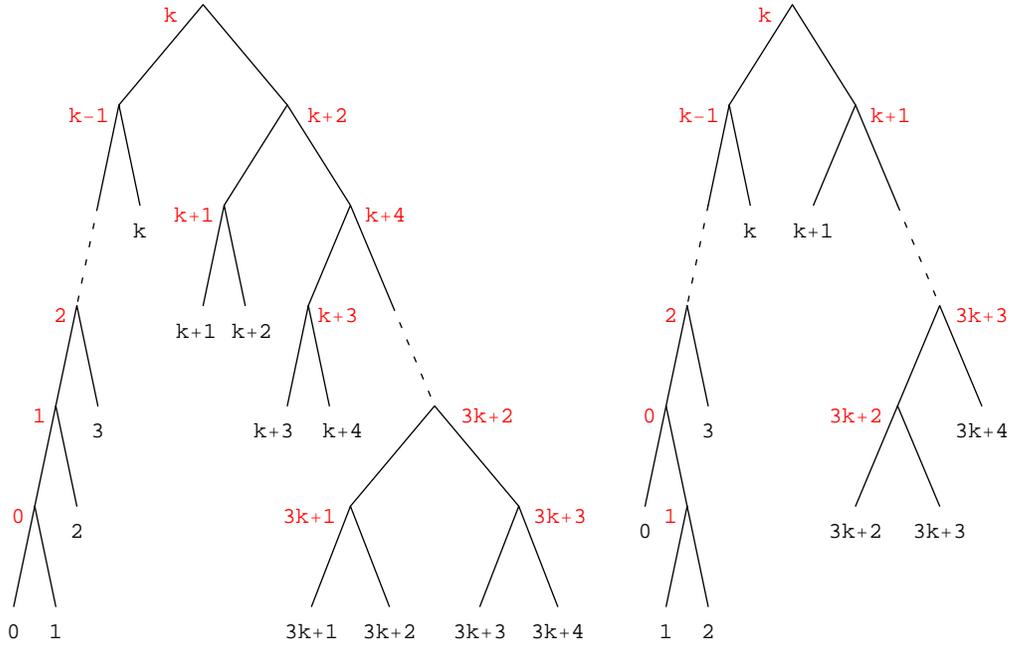}

\caption{The seesaw element $x_0^{k} x_1 x_{3k+3} x_{3k+2}^{-1} x_{3k
 }^{-1} \ldots x_{k+4}^{-1} x_{k+2}^{-1} x_0^{-k-1}$ of swing $k$ in
 $F$.}

\ec
  \label{fig:seesaw}
\end{figure}

Examples of seesaw elements in $F$ are given explicitly in the
following theorem, and illustrated in Figure \ref{fig:seesaw}.
We note that the description given below is in terms of the infinite
generating set $\{x_0, x_1, x_2, \ldots \}$ of $F$ for
convenience but these could be expressed in terms of $x_0$ and
$x_1$ by substituting $x_n= x_1^{x_0^{n-1}}$ for $n \geq 2$.

\begin{thm}[\cite{\seesaw}, Theorem 4.1] \label{thm:seesaw} 
The elements   $$x_0^{k} x_1 x_{3k+3} x_{3k+2}^{-1} x_{3k
}^{-1} \ldots x_{k+4}^{-1} x_{k+2}^{-1} x_0^{-k-1}$$ are seesaw
elements of swing $k$ with respect to the generator $x_0$ in the
standard generating set $\{x_0, x_1\}$.
\end{thm}

We use the seesaw elements defined in Theorem \ref{thm:seesaw} to
find an infinite number of distinct cone types in $F$.

\begin{thm} \label{thm:Finfmanyconetypes} 
Thompson's group $F$ contains infinitely many cone types with respect
to the generating set $\{x_0,x_1\}$.
\end{thm}

\textit{Proof}:    The seesaw
elements given in Theorem \ref{thm:seesaw} are all defined with
respect to the generator $x_0$ of $F$. Let $w$ be a seesaw element of
swing $k$ of the form given above. The possible geodesic continuations
of the word $w x_0^{-l}$ where $l \in [1,k]$ includes $x_0^l$ but not
$x_0^{l+1}$.  Varying $l$, we have produced a finite set of group
elements with distinct cone types.  Varying $k$, the ``swing'' of the
element, we can produce larger and larger finite sets of distinct cone
types, so the set of cone types is unbounded.
\hfill$\Box$

The next theorem follows easily because the relators of the group $F$ in
the finite presentation given above all have even length.

\begin{thm}\label{thm:FnotREG}
The full language of geodesics in $F$ with respect to
$\{x_0,x_1\}$ is not regular.
\end{thm}

\textit{Proof}: 
The finite presentation of $F$ given above has relators of lengths
$10$ and $14$, and Theorem \ref{thm:Finfmanyconetypes} shows that
$F$ has infinitely many cone types.  It then follows from Lemma
\ref{lem:evenlength} that the full language of geodesics is not
regular.
\hfill$\Box$

\section{Groups with infinitely many cone types}
\label{sec:generalcase}

The main theorems of this section show that if a group $G$ with
finite generating set $S$ contains an infinite family of the
seesaw elements of arbitrary swing, as defined in Section \ref{sec:F}, then two results follow:
\bi
\item $G$ has no regular language of geodesics, and 
\item $G$ has infinitely many cone types with respect to $S$.
\ei
These seesaw elements are present
in $F$ with the standard generating set $\{x_0,x_1\}$, in $L_m$
with the wreath product generators, and in a large class of wreath
products as described in \cite{\deadlamp}.

%
%

Seesaw elements of large swing preclude the possibility of there
being a regular language of geodesics, by an argument similar to that used
to prove Theorem \ref{thm:FnotREG}.

\begin{thm} \label{thm:seesawnotreg}
A group $G$ generated by a finite generating set $X$ with seesaw
elements of arbitrary swing with respect to $X$ has no regular
language of geodesics.
\end{thm}

\textit{Proof}: 
Suppose there were a regular language of geodesics for $G$ with
pumping length $p$, and consider the form of the Pumping
Lemma used in Theorem \ref{lem:LnotREG}.   We take $w$ to be a
seesaw element of swing $k$ with respect to a generator $t$ with
$k>p$, and note that any geodesic representative for $w$ must be
written $w = v t^{k}$ or $w=v' t^{-k}$. So a word in one of these
two forms must belong to the regular language.

First suppose that $w = v t^{k}$ belongs to the regular language.
Applying the Pumping Lemma to the suffix $t^{k}$ of $w$, we see
that $v t^{k+n}$ must be in the regular language as well.  Since
$k+n>k$, this path is geodesic only until $v t^{k}$, and after
that, further multiplication by $t$ will decrease word length.
Thus $v t^{k+n}$ cannot be geodesic, contradicting the Pumping
Lemma.

Similarly, if the regular language contains a geodesic
representative of $w$ of the form $v' t^{-k}$, we again apply the
Pumping Lemma to obtain a contradiction. Thus there can be no
regular language of geodesics for $G$.
\hfill$\Box$

Since Thompson's group $F$ contains seesaw elements of arbitrary swing \cite{seesaw}, it follows from Theorem \ref{thm:seesawnotreg} that $F$ has no regular language of geodesics.  Theorems \ref{lem:LnotREG} and \ref{lem:LnotREGauto} also follow from Theorem \ref{thm:seesawnotreg}, since $L_m$ in the wreath product generating set is shown to have seesaw elements of arbitrary swing in \cite{deadlamp}, and $L_2$ is shown to have seesaw elements of arbitrary swing with respect to the automata generating set in \cite{lamptat}.

The following corollary follows immediately from Lemma \ref{lem:coneimpliesreg}.

\begin{cor}\label{cor:conetypes}
Let $G$ be a finitely generated group with finite generating set
$S$. If $G$ has seesaw elements of arbitrary swing then $G$ has
infinitely many cone types with respect to the finite generating
set $S$.
\end{cor}

This corollary provides alternative proofs of Corollaries
\ref{cor:lmconetypes} and \ref{cor:l2conetypes}, and  Theorem
\ref{thm:Finfmanyconetypes}.

Seesaw elements are also found in more general wreath products.

\begin{thm}[\cite{\deadlamp}, Theorem 6.3]
\label{thm:seesawswreath} Let $F$ be a finitely generated group
containing an isometrically embedded copy of $\Z$, and $G$ any
finitely generated non-trivial group. Then
$G \wr F$ contains seesaw elements of arbitrary swing with respect to
at least one generating set.
 \end{thm}

 Combining Theorems  \ref{thm:seesawnotreg} and
 \ref{thm:seesawswreath} with Corollary \ref{cor:conetypes}, we obtain the following corollary.

 \begin{cor}
 Let $F$ be a finitely generated group containing an isometrically
embedded copy of $\Z$, and $G$ any finitely generated non-trivial group. Then $G \wr F$ contains infinitely many cone
types with respect to at least one generating set and there is no
regular language of geodesics for $F$ with respect to that generating
set.
\end{cor}







\end{document}